\RequirePackage{fix-cm}

\documentclass[smallextended]{svjour3}       
\usepackage{amsmath}
\smartqed  
\usepackage{tikz-network}
\usepackage{graphicx}
\usepackage{color}
\usepackage{amssymb} 
\usepackage{framed}
\usepackage{wrapfig}
\usepackage{hyperref}
\usepackage{stmaryrd}
\hypersetup{
	colorlinks=true,
	linkcolor=black,
	urlcolor=red,
}

\newcommand{\F}{\mathsf F}

\newcommand{\Q}{\mathbf Q}
\newcommand{\R}{\mathbf R}

\newcommand{\Z}{\mathbf Z}

\newcommand{\C}{\mathbf C}

\newcommand{\Jimm}{\mathbf J}

\DeclareMathOperator{\con}{con}

\DeclareMathOperator{\num}{num}

\newcommand{\pgl}{ \mathsf{PGL}_2 (\Z)   }

\definecolor{green}{RGB}{117, 165, 50}

\usepackage{float}
\usepackage{forest}
\usepackage{authblk}
\newcommand{\sherh}[1]{\fboxsep=0pt\setlength{\fboxrule}{1pt}
\begin{center}
   \fbox{\colorbox{green}{
         \begin{minipage}[t]{17cm}
            #1
         \end{minipage}
      }
   }
\end{center}}
\newcommand{\sherhh}[1]{\fboxsep=0pt\setlength{\fboxrule}{1pt}
\begin{center}
   \fbox{\colorbox{yellow}{
         \begin{minipage}[t]{17cm}
            #1
         \end{minipage}
      }
   }
\end{center}}

\newcommand{\sherhhh}[1]{\fboxsep=0pt\setlength{\fboxrule}{1pt}
\begin{center}
   \fbox{\colorbox{red}{
         \begin{minipage}[t]{17cm}
            #1
         \end{minipage}
      }
   }
\end{center}}

\renewcommand{\sherh}[1]{}\renewcommand{\sherhh}[1]{}
\renewcommand{\sherhhh}[1]{}

\newcommand{\quant}[2]{\left\llbracket #1 \right\rrbracket_{#2}}
\newcommand{\quan}[1]{\left\llbracket #1 \right\rrbracket}

\begin{document}

\title{Quantizations of  Continued Fractions
}

\titlerunning{{Quantizations of  Continued Fractions}}        

    \author{ A. Muhammed Uluda\u{g}$^a$   \and   Esra \"{U}nal Y{\i}lmaz$^{b*}$  }

\authorrunning{A. Muhammed Uluda\u{g}  \and  Esra \"{U}nal Y{\i}lmaz} 

\institute{A. Muhammed Uluda\u{g}$ ^a $ \email{muhammed.uludag@gmail.com}  \at
Department of Mathematics, Galatasaray University, \c{C}{\i}ra\u{g}an Cad. No. 36, 34357 Be\c{s}ikta\c{s},  \.{I}stanbul, Turkey \and  
 Esra \"{U}nal Y{\i}lmaz$ ^{b*} $  
	\email{esraunal@ibu.edu.tr} (Corresponding author)   \at
Department of Mathematics, Bolu Abant \.{I}zzet Baysal University,  14100 G\"{o}lk\"{o}y,  Bolu, Turkey \\
}

\date{Received: date / Accepted: date}

\maketitle

\begin{abstract}
 We introduce a four-parameter deformation of continued fractions, which we call $ U $-deformation. We study some particular cases and compare them with the q-deformation of continued fractions introduce recently by  Morier-Genoud and Ovsienko. 
\end{abstract}

\section{Introduction}\label{intro}
Consider the \emph{numerator  function} $\num:\Q^+\to \Z^+$ defined as
\begin{equation}
\num: x=(p,q)\in \Q^+\to \num(x)=p \in \Z^+,
\end{equation}
where $p, q>0,$ $ \gcd(p,q)=1$. 
It is easy to verify that $\num$ satisfies the following system of equations 
for a function $f:\Q^+\to \C$:
\begin{eqnarray}\label{numerator}
f(1+x)=f(x)+f(1/x), \label{nmr111}\\
f\left(\frac{x}{1+x}\right)=f(x), \label{nmr222}
\end{eqnarray}
and $\num$ is the unique solution of this system satisfying $f(1)=1$.
The current paper is devoted to studying the following generalization of this system:
\begin{eqnarray}\label{numerator2}
f(1+x)=pf(x)+qf(1/x), \label{eq:M1}\\
f\left(\frac{x}{1+x}\right)=rf(x)+sf(1/x), \label{eq:M2}
\end{eqnarray}
where $f:\Q^+\to \C$ and $p,q,r,s\in\C$. Its solution is linear in $\C$; so there is no danger in assuming $f(1)=1$ from now on. With this assumption, one has
$$
f(1)=1, \quad f(2)=p+q, \quad f(1/2)=r+s.
$$
Consider the monoid $\mathcal M$ generated by the transformations $x\to 1+x$ and $x\to x/(1+x)$. By using the simple continued fraction expansion $[n_0,n_1, \cdots, n_k]$ of an $x\in \Q^+$, 
it is easy to show that there exists a unique $M\in \mathcal M$ with $M(1)=x$.
Set $\ell(x):=n_0+n_1+ \cdots+n_k$.
Since $\ell(1+x)=\ell(x/(1+x))=\ell(x)+1$, 
equations (\ref{eq:M1}-\ref{eq:M2}) expresses $f$ in terms of $f$-values of rational numbers with 
smaller $\ell$-value. Hence (\ref{eq:M1}-\ref{eq:M2}) are consistent and $f(x)$ can be computed in terms of $f(1)$. Therefore a solution always exists and is unique; if we further impose 
$f(1)=c$. 
Furthermore, $f$ is an integral-valued function  if $p,q,r,s\in \Z$ and
positive valued  if $p,q,r,s\in \Z^+$.

Clearly, the system (\ref{eq:M1}-\ref{eq:M2})
reduces to (\ref{nmr111}-\ref{nmr222}) when $(p,q,r,s)=(1,1,1,0)$. Moreover,~(\ref{eq:M1}-\ref{eq:M2}) implies
\begin{equation}\label{recc}
f(2+x)=pf(1+x)+qf(1/(1+x))=pf(1+x)+qrf(1/x)+qsf(x).
\end{equation}
When $r=0$, this gives 
$$
f(2+x)=pf(1+x)+qs f(x),
$$
i.e~(\ref{eq:M1}-\ref{eq:M2}) can be viewed as an extension of the Fibonacci recursion
to rational indices.

We can rewrite (\ref{eq:M1}-\ref{eq:M2}) as
\begin{eqnarray}\label{numeratorm}
\left(
f(1+x),
f\left({x}/{(1+x)}\right)
\right)
=
U
\left(\begin{matrix}
f(x)\\
f\left({1}/{x}\right)
\end{matrix}\right), \quad 
U=\left(\begin{matrix}
p&q\\
r&s
\end{matrix}\right)
\end{eqnarray}
where $f:\Q^+\to \C$ and $ f(1)=1$.
Denote its solution by $f_U$. 
In this paper we suggest the function $\quant{x}{U}:=\frac{f_U(x)}{f_U(1/x)}$ as a ``$ U$-deformation" of the real number $ x $ and compare some special cases with the $ q $-deformation of continued fractions introduced by Morier-Genoud and Ovsienko \cite{Valentin}. 

In particular, when  $U={\tiny \left(\begin{matrix}
1\!\!&1\\
1\!\!&0
\end{matrix}\right)}$,
we get the system~(\ref{nmr111}-\ref{nmr222}). Thus, in this case,
$f_U=\num$. Another particular case is $U={\tiny \left(\begin{matrix}
1\!\!&1\\
0\!\!&1
\end{matrix}\right)}$. In this case the system~(\ref{eq:M1}-\ref{eq:M2}) becomes
\begin{eqnarray}\label{connumerator2}
f(1+x)=f(x)+f(1/x), \label{eq:conM1}\\
f\left({x}/{(1+x)}\right)=f(1/x). \label{eq:conM2}
\end{eqnarray}
Its solution $\con:=f_{U}$
was called the {\it conumerator} and the function $\F(x):=\con(1/x)$ was called the {\it codenominator} in~\cite{conumerator}. The latter function solves therefore the system
\begin{eqnarray}\label{deconnumerator2}
\F(1/(1+x))=\F(x)+\F(1/x), \label{eq:deconM1}\\
\F\left(1+1/x\right)=\F(x). \label{eq:deconM2}
\end{eqnarray}
One has
$$
\F(2+x)=
\F(1+\frac{1}{1/(x+1)})\stackrel{(\ref{eq:deconM2})}{=}
\F(1/(1+x))\stackrel{(\ref{eq:deconM1}-\ref{eq:deconM2})}{=} \F(1+x)+\F(x),
$$
i.e. $\F$ satisfies the Fibonacci recursion as well. Indeed, one can show that
$\F(n)=F_n$, where $F_n$ is the $n$th Fibonacci number defined by $F_1=F_2=1$ and $F_{n+2}=F_{n+1}+F_n$. Hence, $\F$ gives an extension of the Fibonacci sequence to $\Q^+$ as an integral-valued map, i.e. $\F(4/9)=11$ is the $4/9$th Fibonacci number.
Moreover, the function
$$
\Jimm:\Q^+\to \Q^+,\quad \Jimm(x):=\frac{\con(x)}{\con(1/x)}=\frac{\F(1/x)}{\F(x)}
$$
is an involution, induced by Dyer's outer automorphism of the group $\pgl$. It can be continuously extended to $\R\setminus \Q$, and satisfies the functional equations~\cite{conumerator}
$$
\Jimm(1/x)=1/\Jimm(x), \quad \Jimm(1-x)=1-\Jimm(x), \quad \Jimm(-x)=-1/\Jimm(x)
$$
on a certain subset of $\R\setminus \Q$.
If  $x=[n_0,n_1, \dots, n_k ]$ is given as a simple continued fraction, then
\begin{equation}\label{acord}
		\Jimm(x)= [1_{n_0-1},2, 1_{n_1-2},2,1_{n_2-2},2, 
		\dots 
		2,1_{n_{k-1}-2},2,1_{n_k-1}],
		 \end{equation}
 where  $1_{k}$ denotes the sequence $1,1,\dots , 1$ of length $k$.
(It is understood that $[ \dots , n, 1_0, m, \dots ]:= [ \dots, n,m, \dots ]$ and
	$[ \dots , n, 1_{-1}, m, \dots ]:=  [ \dots , n+m-1, \dots ].$)

\section{Quantization}
Consider the function
$$
\quant{x}{U}:=\frac{f_U(x)}{f_U(1/x)} \implies \quant{\frac{1}{x}}{U}=\frac{1}{\quant{x}{U}}
$$
Then 
$$
\quant{1}{U}=1, \quad \quant{2}{U}=\frac{p+q}{r+s}, \quad \quant{\frac12}{U}=\frac{r+s}{p+q}.
$$
In particular, one has
$$
\quant{x}{{\tiny \left(\begin{matrix}
1\!\!&1\\
0\!\!&1
\end{matrix}\right)}}=x, \quad 
\quant{x}{{\tiny \left(\begin{matrix}
1\!\!&1\\
1\!\!&0
\end{matrix}\right)}}=\Jimm(x),
 \quad 
\quant{\Jimm(x)}{{\tiny \left(\begin{matrix}
1\!\!&1\\
1\!\!&0
\end{matrix}\right)}}=x.
$$
Our aim is to study the limit $\lim_{x\to y}\quant{x}{U}$ when $x\in \Q$, $y\in \R$ and to propose this limit as a ``{quantization of $y$ with parameters $p,q,r,s$}'' (or as ``{$U$-deformation}'') of $y$. 

From now on, we shall suppress $U$ from $\quant{x}{U}$ and write simply $\quan{x}$ 
unless it is necessary to indicate $U$. Whence
$$
\quan{1+x}=\frac{pf_U(x)+qf_U(1/x)}{rf_U(1/x)+sf_U(x)}
=\frac{p\quan{x}+q}{r+s\quan{x}} \implies
$$
\begin{equation}\label{funk}
\quan{1+x}=
\begin{cases}
\frac{p}{s}+\frac{qs-rp}{s(s\quan{x}+r)} & (s\neq 0)\\
\frac{p}{r}\quan{x}+\frac{q}{r} & (s=0).
\end{cases}
\end{equation}
Set $\Delta=\Delta(U):=qs-rp$. Quite naturally, assume $\Delta\neq 0$. 

\section{\bf The case $s=0$.} Set $P:=p/r$ and $Q:=q/r$.  
Then by iterating~(\ref{funk}) we get
\begin{equation}\label{funkt}
\quan{x+n}=P^n \quan{x}+\frac{1-P^n}{1-P}Q.
\end{equation}
Let 
$$
x=n_0+\cfrac{1}{n_1+\cfrac{1}{\dots}}=:[n_0, n_1, \dots, n_k]
$$
be a finite or infinite simple continued fraction.
Then by using~(\ref{funkt}) and the equation $\quan{1/x}=1/\quan{x}$, we can express 
$\quan{x}$ as a continued fraction as
\begin{eqnarray}\label{numerator34}
\quan{n_0,n_1,n_2,\dots}=\frac{1-P^{n_0}}{1-P}Q+
\cfrac{P^{n_0}}
{\frac{1-P^{n_1}}{1-P}Q +
\cfrac{
P^{n_1}}{\frac{1-P^{n_2}}{1-P}Q +\dots}}.
\end{eqnarray}

when $P\neq 1$ and as
\begin{eqnarray}\label{numerator124}
\quan{n_0,n_1,n_2,\dots}=
n_0Q+\cfrac{1}{n_1Q +\cfrac{1}{n_2Q +\dots}}
=[n_0Q,n_1Q,n_2Q,\dots]
\end{eqnarray}
when $P=1$.
In particular, we have
\begin{eqnarray}\label{numerator44}
\quan{1,1,1,\dots}=Q+
\cfrac{P}
{Q +
\cfrac{P}{Q +\cfrac{P}{\dots}}}=\frac{Q+\sqrt{Q^2+4P}}{2}
\end{eqnarray}
For the moment we don't worry about the convergence of these continued fractions and the topology of this convergence, which can be taken to be Archimedean or non-Archimedean.
 When $ Q=1 $ we obtain
\begin{eqnarray} \label{confr}
\quan{n_0,n_1,n_2,\dots}=\frac{1-P^{n_0}}{1-P}+
\cfrac{P^{n_0}}
{\frac{1-P^{n_1}}{1-P} +
\cfrac{
P^{n_1}}{\frac{1-P^{n_2}}{1-P} +\dots}}.
\end{eqnarray}
Compare this with the q-deformation of Ovsienko and Morier-Genoud \cite{Valentin} (with their notation)
\begin{equation}
\label{qa}
[a_{1}, \ldots, a_{2m}]_{q}:=
[a_1]_{q} + \cfrac{q^{a_{1}}}{[a_2]_{q^{-1}} 
          + \cfrac{q^{-a_{2}}}{[a_{3}]_{q} 
          +\cfrac{q^{a_{3}}}{[a_{4}]_{q^{-1}}
          + \cfrac{q^{-a_{4}}}{
        \cfrac{\ddots}{[a_{2m-1}]_q+\cfrac{q^{a_{2m-1}}}{[a_{2m}]_{q^{-1}}}}}
          } }} 
\end{equation}
where the standard notation 
 $[a]_q=1+q+q^2+\cdots+q^{a-1}$ is used for the $q$-deformation of integer $a$.
\subsection{The case $(p,q,r,s)=(p,1,1,0)$}
In this case, the system (\ref{eq:M1}-\ref{eq:M2}) becomes
\begin{eqnarray*}\label{numerator22}
f(1+x)=pf(x)+f(1/x), \label{eq:M12}\\
f\left(\frac{x}{1+x}\right)=f(x), \label{eq:M22}
\end{eqnarray*}
and \eqref{confr} becomes (with the initial condition  $ f(1)=1 $)
\begin{eqnarray} \label{confr-p}
\quan{n_0,n_1,n_2,\dots}=\frac{1-p^{n_0}}{1-p}+
\cfrac{p^{n_0}}
{\frac{1-p^{n_1}}{1-p} +
\cfrac{
p^{n_1}}{\frac{1-p^{n_2}}{1-p} +\dots}}.
\end{eqnarray}

For example,
$$
\begin{array}{|l|l|}
\hline
x&f_x(p)\\\hline
1/3 \quad  &1\\\hline
1/4 \quad &1\\\hline
3/4&p^{2}+p+1\\\hline
4/5&p^3+p^2+p+1\\\hline
5/6&p^4+p^3+p^2+p+1\\\hline
17/31&p^5+3p^4+3p^3+5p^2+4p+1 \\\hline
17/2&2p^{8}+2p^{7}+2p^{6}+2p^{5}+2p^{4}+2p^{3}+2p^{2}+2p+1  \\\hline
19/34&p^{5}+2p^{4}+5p^{3}+6p^{2}+4p+1 \\\hline
29/13&p^{6}+4p^{5}+7p^{4}+7p^{3}+6p^{2}+3p+1 \\\hline
\end{array}
$$
{\bf Observation 1 (Unimodality).} The coefficients $ a_i $ of $ f_x(p)=\sum a_ip^i $ first increase as $ i $ grows, then they decrease as $ i $ grows. 

\medskip 

As for the value of $\quan{x}$ for $(p,q,r,s)=(p,1,1,0)$, first note that when $ n $  is a positive integer, its $ U $ deformation coincides with the usual definition of $ p$-deformation:
$$\quan{n}=\frac{1-p^n}{1-p}. $$
{\bf Observation 2} For general $ x\in \Q^+, $ experiments indicate that the Taylor expansion of $\quan{x}$  around the origin always has integral coefficients.
For example
$$
\quan{\frac{7}{5}}=\frac{3p^2+3p+1}{2p^2+2p+1}
$$
$$=1 + p- p^2 + 2p^4 - 4p^5 + 4p^6 - 8p^8 + 16p^9 - 16p^{10} + 32p^{12} - 64p^{13} + 64p^{14} + O(p^{15})$$

(Compare this with the quantization of \cite{Valentin} (with their notation))
$$
\left[\frac{7}{5}\right]_q=\frac{1+q+2q^2+2q^3+q^4}{1+q+2q^2+q^3}
$$ 
\begin{align*}
=1+q^3-2q^5+q^6+3q^7-3q^8-4q^9+7q^{10}+4q^{11}-14q^{12}.
\end{align*} Note that the polynomials $ f_x(p) $ has positive integral coefficients. This is obvious for us by definition. On the other hand, in \cite{Valentin} this fact was called {\it total positivity}, it is not obvious and had to be proven.

 Let us give another example
 $$
\quan{\frac{19}{31}}=\frac{5p^3+8p^2+5p+1}{2p^4+10p^3+12p^2+6p+1}
$$
\begin{align*}
=1-p+2p^{2}-5p^{3}+14p^{4}-42p^{5}+130p^{6}-406p^{7}+1268p^{8}-3952 p^{9}
\\+12296p^{10} 
-38220p^{11}+118752p^{12}-368928p^{13} 
+1146152p^{14}+\mathrm{O}\! \left(p^{15}\right)
\end{align*}
(Compare this with the quantization of \cite{Valentin} (with their notation))
$$
\left[\frac{19}{31}\right]_q=\frac{2q^6+4q^5+5q^4+4q^3+3q^2+q}{2q^6+4q^5+7q^4+7q^3+6q^2+4q+1}
$$
\begin{align*}
=q-4 q^{2}+14 q^{3}-34 q^{4}+77 q^{5}-173 q^{6}+384 q^{7}-847 q^{8}+1864 q^{9} 
\\-4091 q^{10}+8959 q^{11}-19599 q^{12}+42851 q^{13}-93648 q^{14}+\mathrm{O}\! \left(q^{15}\right).
\end{align*}
{\bf{Observation 3}} Taylor series of $ f_x(p) $ has alternating coefficients.

 The simplest example of an infinite continued fraction is the expansion of the golden ratio $$\frac{1+\sqrt{5}}{2}=[1,1,1,1,\dots].  $$ According to \eqref{confr-p}
\begin{equation}\label{conf1}
\quan{1,1,1,1,\dots}=1+
\cfrac{p}
{1+
\cfrac{
p}{1 +\dots}}
\end{equation} 
Taylor series of \eqref{conf1} starts as follows:
\begin{align*}
=&1+p-p^{2}+2p^{3}-5p^{4}+14p^{5}-42p^{6}+132p^{7}-429p^{8}+1430p^{9}-4861 p^{10} \\
+&16778p^{11}-58598p^{12}+206516p^{13}-732825p^{14}+2613834 p^{15}-9358677p^{16} \\
+&33602822p^{17}-120902914p^{18}+435668420p^{19}+\mathrm{O}\! \left(z^{20}\right). 
\end{align*}
We were able to identify the coefficients $a_n$ appearing in this series as the so-called Catalan numbers (see sequence A000108 of OEIS), but with alternating signs (see \cite{mythesis}).
Note that the coefficients of the Taylor expansion of the quantization of \cite{Valentin} (with their notation)
\begin{align*}
[1,1,1,\dots]_q=&1+q^{2}-q^{3}+2q^{4}-4q^{5}+8q^{6}-17q^{7}+37q^{8}-82q^{9}+185q^{10} \\
-&423q^{11}+978q^{12}-2283q^{13}
+5373q^{14}-12735q^{15}+30372q^{16} \\
-&72832q^{17}+175502q^{18}-424748q^{19}+1032004q^{20}\dots
\end{align*}
are the so-called Generalized Catalan numbers (see sequence A004148 of OEIS), but with alternating signs. 

Observation 2 is true in general: 
\begin{theorem} \label{integerth}
The Taylor expansion   $ \quant{x}{\tiny \left(\begin{matrix}
p\!\!&1\\
1\!\!&0
\end{matrix}\right)}=\sum_{n=0}^{\infty}a_np^n  $ has integral coefficients.
\end{theorem}
\begin{proof}
If we consider  \eqref{confr-p}, we can write   
\begin{align} \label{intcoef11}
\quan{x_k}&=\quan{{n_0,n_1,\dots,n_k}}=\frac{1-p^{n_0}}{1-p}
          + \cfrac{p^{n_{0}}}{\frac{1-p^{n_1}}{1-p} 
          + \cfrac{p^{n_{1}}}{
        \cfrac{\ddots}{\frac{1-p^{n_{k-1}}}{1-p}+\cfrac{p^{n_{k-1}}}{\frac{1-p^{n_{k}}}{1-p}}}}}.
\end{align}
Obviously $ \frac{1}{1+pf(p)} $ has integral Taylor coefficients if $f$ does.
If we apply with this rule inductively we get the desired result. For example, for the induction basis we have 
$$
\frac{1-p^{n_k}}{1-p}=1+p(1+p+p^2+\dots+p^{n_{k-2}})
$$ 
with $ f(p)= 1+p+p^2+\dots+p^{n_{k-2}}$.
\end{proof}
\begin{proposition}(Stabilization phenomenon).
Let $ x\geq 1 $  be an irrational real number. The Taylor expansions at $p=0  $ of two
consecutive {\tiny $\left(\begin{matrix}
p\!\!&1\\
1\!\!&0
\end{matrix}\right)$}-deformed convergents of the continued fraction of $ x $, namely of $ x_{k-1}=[n_0,\dots,n_{k-1}]$ and
$ x_{k}=[n_0,\dots,n_k]$, have the first $ n_0+n_1+\dots+n_{k}-1 $ terms identical.
\end{proposition}
\begin{proof}
Let  $ x\geq 1 $ be irrational real number. Let $ x_{k}=[n_0,...,n_k]$. Consider two sequences $ (R_k)_{k\geq0} $ and $ (S_k)_{k\geq0} $ of polynomials in $ p$ defined by the following recursions
\begin{align*}
R_{k+1}=\quan{n_k}R_k+p^{n_{k-1}}R_{k-1} \\
S_{k+1}=\quan{n_k}S_k+p^{n_{k-1}}S_{k-1}
\end{align*}
with the initial conditions $ (R_0,R_1)=(1,\quan{n_0}) $ and $ (S_0,S_1)=(0,1) $ where polynomials $ R $ and $ S $ represent quantization polynomials $ f(x) $ and $ f(1/x)$ in $ p. $ Hence 
$$\quan{x_{k-1}}=\frac{R_{k-1}}{S_{k-1}},  \quad \quan{x_k}=\frac{R_k}{S_k}.  $$ 
So that
\begin{equation}
 \quan{x_k}-\quan{x_{k-1}}=\frac{R_k}{S_k}-\frac{R_{k-1}}{S_{k-1}}=\frac{R_kS_{k-1}-S_kR_{k-1}}{S_kS_{k-1}}.
\end{equation}
The polynomial in the numerator of the right hand side is 
\begin{align*}
R_kS_{k-1}-S_kR_{k-1}=&(\quan{n_{k-1}}R_{k-1}+q^{n_{k-2}}R_{k-2})(\quan{n_{k-2}}S_{k-2}+q^{n_{k-3}}S_{k-3}) \\
-& (\quan{n_{k-1}}S_{k-1}+q^{n_{k-2}}S_{k-2})(\quan{n_{k-2}}R_{k-2}+q^{n_{k-3}}R_{k-3}) \\
=&(-1)^kp^{n_0+...+n_{k-2}}.
\end{align*}
Both polynomials $ S_k $ and $ S_{k-1} $ start with zero order term $ 1 $ so is the series $ 1/S_kS_{k-1}. $ Then $$ \quan{x_k}-\quan{x_{k-1}}=(-1)^kp^{n_0+...+n_{k-2}}+\mathcal{O}({p^{n_0+...+n_{k-2}+1}}). $$ 
 This shows that the coefficients stabilize.

\end{proof}

\subsection{Quantization of $e$ and $\pi$}
In this section we write down the first terms of the quantization of two notable examples of transcendental irrational numbers, $e$ and $\pi$
The continued fraction expansion of the Euler constant is given by the following famous regular pattern $e=[2,1,2,1,1,4,1,1,6,1,1,8,1,1,10,\dots].  $
The Taylor series of $ \quan{e} $ starts as follows: 
\begin{align*}
\quan{e}=&1+p+p^{2}+p^{3}+p^{4}+p^{5}-p^{6}+2p^{8}-3p^{9}+p^{10}+2p^{11}-p^{12} \\
-&2 p^{13}-7p^{14}+32p^{15}-36p^{16}-25p^{17}+86p^{18}+40p^{19}-351p^{20} \\
+&297 p^{21}+802p^{22}-2528p^{23}+3422p^{24}-4226p^{25}+10375p^{26} \\
-&24836 p^{27}+32704p^{28}-408p^{29}-96389p^{30}+266726p^{31} \\
-&603941p^{32}+1426546p^{33}-3149939p^{34}+5686426p^{35}-7795364p^{36} \\
+&6880187p^{37}+2141014 p^{38}-32473209p^{39}+\mathrm{O}\! \left(z^{40}\right)
\end{align*}
We observe that the coefficients of $ p^{17+2k} $ where $k\geq 0 $, turn out to be smaller than those of their neighbours (compare \cite{Valentin}). 

The continued fraction expansion of $\pi  $ 
starts as follows: 
$$\pi=[3,7,15,1,292,1,1,1,2,1,3,1,14,2,1,1,2,2,2,2,1,84,\dots]. $$
The Taylor series of $ \quan{\pi} $ starts as follows: 
\begin{align*}
\quan{\pi}=&1+p+{p}^{2}+{p}^{3}-{p}^{4}+2\,{p}^{11}-3\,{p}^{12}+{p}^{13}+4\,{p}^{
19}-8\,{p}^{20}+5\,{p}^{21}-{p}^{22} \\
-&2\,{p}^{26}+17\,{p}^{27}-40\,{p}^
{28}+52\,{p}^{29}-62\,{p}^{30}+90\,{p}^{31}-144\,{p}^{32}+233\,{p}^{33
} \\
-&385\,{p}^{34}+666\,{p}^{35}-1133\,{p}^{36}+1829\,{p}^{37}-2904\,{p}^
{38}+4656\,{p}^{39}+O \left( {p}^{40} \right). 
\end{align*}
The coefficients of above series grow very slowly in contrast with the other examples we considered so far and we observe that the coefficients increase from coefficients of $ p^{27} $ regardless of signs.

\section{The case $r=0$.}
\subsection{The case $(p,q,r,s)=(p,1,0,1)$}
Here~(\ref{recc}) become $f(2+x)=pf(1+x)+qsf(x)$ with 
$f(1)=1$, $f(2)=p+q$. When $q=1$, $f(n)$ is a kind of Fibonacci polynomial in $p$, with initial values $f(1)=1$ and $f(2)=1+p$ (whereas for the official Fibonacci polynomial one has $f(1)=1$, $f(2)=p$).  Recall that the Fibonacci polynomials are a polynomial sequence which can be considered as a generalization of the Fibonacci numbers. These Fibonacci polynomials are defined by a recurrence relation:
$$F_n(p)=pF_{n-1}(p)+F_{n-2}(p)  $$ with initial values $ F_0(p)=0, \quad F_1(p)=1 \implies F_2(p)=1+p.$ Hence, $ f$ equals the Fibonacci polynomial $ F_x(p) $ for integral $ x. $
For rational values of $x$, $f$ defines an extension of Fibonacci polynomial with rational index $ x $.
For example,
$$
\begin{array}{|l|l|}
\hline
x&f_x(p)\\\hline
F_n/F_{n+1} \quad (n=1,2,\dots) &1\\\hline
F_{n+1}/F_n \quad (n=1,2,\dots)&(n-1)p+1\\\hline
3/4&p+1\\\hline
4/3&2p^2+2p+1\\\hline
4/5&p^2+p+1\\\hline
17/31&3p^3+p^2+4p+1\\\hline
17/2&2p^8+p^7+13p^6+6p^5+25p^4+10p^3+14p^2+4p+1 \\\hline
\end{array}
$$
{\bf Observation (Anti-unimodality).} The coefficients $ a_i $ of $ f_p=\sum a_ip^i $ satisfies   $ a_i \geq a_{i+1} $  if $ i $ is even and  $ a_i \leq a_{i+1} $  if $ i $ is odd.

As for the value of $\quan{x}$ for $(p,q,r,s)=(p,1,0,1)$, experiments indicate that the Taylor expansion of $\quan{x}$ around 0 has integral coefficients. For example
\begin{align*}
\quan{\frac{17}{2}}=&\frac{2p^8+p^7+13p^6+6p^5+25p^4+10p^3+14p^2+4p+1}{2p^7+p^6+11p^5+5p^4+16p^3+6p^2+5p+1} \\
=&1-p+13p^2-65p^3+283p^4-1233p^5+5465z^6-24273p^7+107594p^8+\cdots
\end{align*}
Indeed, we have the following theorem:
\begin{theorem}
The Taylor expansion   $ \quant{x}{\tiny \left(\begin{matrix}
p\!\!&1\\
0\!\!&1
\end{matrix}\right)}=\sum_{n=0}^{\infty}a_np^n  $ has integral coefficients.
\end{theorem}
\begin{proof}
In this case $ \quan{1+x} $ transforms into
\begin{equation} \label{func2}
\quan{1+x}=
\begin{cases}
\frac{p}{s}+\frac{qs-rp}{s(s\quan{x}+r)} & (r\neq 0)\\
\frac{p}{s}+\frac{q}{s\quan{x}} & (r=0).
\end{cases}
\end{equation}
Then from \eqref{func2} we get
\begin{eqnarray*}\label{conumerator} 
\quan{n_0,n_1,n_2,\dots,n_k}=p+\cfrac{1}{\quan{n_0-1,n_1,\dots ,n_{k-1},n_k}} \\
=p+\cfrac{1}{p+\cfrac{1}{\quan{n_0-2,n_1,n_2,\dots ,n_k}}}
=\dots =p+\cfrac{1}{p+\cfrac{1}{\ddots{2p+\cfrac{1}{\quan{n_1-1,n_2,\dots ,n_k}}}}} \\
=\dots =p+\cfrac{1}{p+\cfrac{1}{\ddots{2p+\cfrac{1}{\quan{n_2,\dots ,n_k}}}}}   =\dots =p+\cfrac{1}{p+\cfrac{1}{\ddots{2p+\cfrac{1}{p+\cfrac{1}{\quan{n_{k-1}-1,n_k}}}}}} \\
=\dots  =p+\cfrac{1}{p+\cfrac{1}{\ddots{2p+\cfrac{1}{p+\quan{n_k}}}}}=\dots = p+\cfrac{1}{p+\cfrac{1}{\ddots{2p+\cfrac{1}{\ddots+\cfrac{1}{p+1}}}}}
\end{eqnarray*} 
In the above continued fraction expansion the last term is $ np+\frac{1}{p+\quan{1}},$ where $ n $ is some natural number. The same argument as in the proof of in Theorem \ref{integerth} we conclude that Taylor expansion has integral coefficients. 
\end{proof}
Note that in the above continued fraction expansion of the total number of $ p $'s is $ n_0+n_1+n_2+\dots+n_k-1. $

\vfill

\newpage
\renewcommand{\arraystretch}{1}

\end{document}